\documentclass[12pt,a4paper]{article}

\usepackage{latexsym}
\usepackage{amssymb}

\usepackage{amsmath}

\newfont{\bl}{msbm10 scaled \magstep1}
\newfont{\bs}{msbm9}
\newfont{\bss}{msbm7}

\numberwithin{equation}{section}

\newcommand{\R}{\mbox{{\bl R}}}
\newcommand{\RR}{\mbox{{\bs R}}}
\newcommand{\RRR}{\mbox{{\bss R}}}
\newcommand{\N}{\mbox{{\bl N}}}

\newcommand{\C}{\mbox{{\bl C}}}
\newcommand{\eps}{\varepsilon}
\newcommand{\be}{\begin{equation}}
\newcommand{\ee}{\end{equation}}
\newcommand{\e}{{\cal E}(\Omega)}
\newcommand{\ds}{{\cal D}'(\Omega)}
\newcommand{\re}{^\rho{\cal E}(\Omega)}
\newcommand{\os}{\tilde{\Omega}}
\newcommand{\ro}{^\rho \Omega}
\newcommand{\ros}{^\rho \tilde{\Omega}}

\def\texts#1{\text{\quad #1 \quad}}

\catcode`\@=11
\def\ps@plain{\let\@mkboth\@gobbletwo
     \def\@oddhead{}\def\@oddfoot{\ifodd\count0 \rm\hfil\thepage
     \else\rm\thepage\hfil\fi}
        \def\@evenhead{}\def\@evenfoot{\rm\thepage\hfil}}
\catcode`\@=13
\advance\oddsidemargin by 4mm
\newcommand{\address}[3]{{\sc #1}\\ #2, E-mail: #3}

\begin{document}

{\LARGE  \bf 
Pointwise Values and Fundamental Theorem in the Algebra of Asymptotic 
Functions}
\par
\vspace{1cm}
%
{\sc \large  Todor. D. Todorov}
%
%
\par
\vspace{1.6cm}
%
%
%
{\bf Introduction:} The algebra of asymptotic functions
$^\rho {\cal E}(\Omega)$ on an  open
set $\Omega \subset \R^d$ was introduced by M. Oberguggenberger [8]
and the author
of  this  paper [14] in the framework of A. Robinson's  nonstandard
analysis.  It  can be described as a differential  associative  and
commutative  ring (of generalized functions) which  is  an  algebra
over the field of A. Robinson's asymptotic numbers $^\rho \C$ (A. Robinson
[13] and A. H. Lightstone and A. Robinson [6]). Moreover,
$^\rho {\cal E}((\Omega)$ is supplied with  the chain of imbeddings:
\begin{equation}
    {\cal E}(\Omega) \subset {\cal D}'(\Omega)
    \subset\,^\rho{\cal E}(\Omega)\,,
\end{equation}
where ${\cal E}(\Omega)$ denotes the differential ring of the $C^\infty$-functions
(complex valued) on $\Omega, {\cal D}(\Omega)$ denotes the differential ring of the
functions in ${\cal E}(\Omega)$ with compact support in $\Omega$ and $\ds$ denotes the
differential linear space of Schwartz distributions on $\Omega$. Here $\e\subset\ds$
is the usual imbedding in the sense of distribution theory
(V. Vladimirov [16]). The imbedding $\ds\subset\, \re$, constructed in
(M. Oberguggenberger and T. Todorov [10]), preserves all linear
operations, including partial differentiation of any order, and the
pairing between $\ds$ and ${\cal D}(\Omega)$. Finally, the imbedding $\e
\subset\, \re$
preserves all differential ring operations. In addition, if $T^d$
denotes the usual topology on $\R^d$, then the family $\{\re\}_{\Omega\in T^d}$ is a
sheaf of differential rings and the imbeddings in (0.1) are sheaf
preserving. That means, in particular, that the restriction $F|\Omega'$ is
a well defined element in $^\rho{\cal E}(\Omega')$ for any $F\in\, \re$ and 
any open
$\Omega'\subseteq\Omega$. On these grounds we consider the asymptotic functions in $\re$
as ``generalized functions on $\Omega$": they are ``generalized" functions
(rather than ``classical") because they are not mappings from $\Omega$ into
$\C$. On the other hand, they are still ``functions" (although,
``generalized" ones) : {\bf a)} because of the imbeddings (0.1) and
{\bf b)} because of the sheaf properties mentioned above. As a result,
$\re$, supplied with the imbeddings (0.1), offers a solution to the
problem for multiplication of Schwartz distributions - similar
to the solution to the same problem presented not long ago by J. F.
Colombeau (and other authors) in the framework of standard analysis
(J. F. Colombeau [1]-[2]). Our interest in the algebra $\re$ is due
to its importance for the theory of partial differential equations
- linear (T. Todorov [15]) and nonlinear (M. Oberguggenberger [9]).
\vskip 0,2 true cm

   The purpose of this paper is to show that $\re$ is isomorphic to
a subring of the ring $(^\rho\C)^{\ros}$ of the pointwise
functions from $\ros$ into $^\rho\C$. Here $^\rho\R$ and $^\rho\C$ are the
fields of the $A$. Robinson's asymptotic numbers (A. Robinson [13]),
real and complex, respectively, $^\rho\R^d$ is the corresponding vector
space, and $\ros = \Omega + {\cal I}(^\rho\R^d)$ denotes the set of the
points in $^\rho\R^d$ of the form $x+h$, where $x\in\Omega$ and $h$ is an
infinitesimal point in $^\rho\R^d$. As a result, every asymptotic function
in $\re$ can be identified with a unique pointwise function defined
on $\ros$ taking values in $^\rho\C$. In the particular case of
classical smooth functions our evaluation coincides with the usual
evaluation. Moreover, since $\re$ contains a copy of the  space of
Schwartz distributions, it follows that the distributions are also
pointwise functions taking values in $^\rho\C$. We also prove the
Fundamental Theorem of Integration Calculus in $\re$ and show that
the set of the scalars of $\re$ coincides exactly with the field
$^\rho\C$ and the set of the real scalars of $\re$ coincides exactly with
the field $^\rho\R$.
\vskip 0,2 true cm

   Recall that if $L$ is a differential linear space, then the set
of the scalars ${\cal S}$ in $L$ consists of the elements $F$ in $L$ with
zero partial first derivatives $\partial^\alpha F = 0, |\alpha| = 1$, in symbols,
${\cal S}=\{f\in L:\nabla F=0\}$. If, in addition, $S$ is supplied with
complex conjugation, then the set of the real scalars in $L$ is
given by Re$({\cal S})=\{\pm|F|:F\in {\cal S}\}$.
\vskip 0,2 true cm

    The theory of asymptotic functions is similar to, but somewhat
different from, J. F. Colombeau's nonlinear theory of generalized
functions, where we also have a chain of imbeddings $\e\subset\ds\subset{\cal G}(\Omega)$
with similar properties. Here ${\cal G}(\Omega)$ denotes a typical algebra
of J. F. Colombeau's generalized functions (J. F. Colombeau [1]-[2]).
The set of the scalars $\bar{\C}$, and the real
scalars $\bar{\R}$ of ${\cal G}(\Omega)$ are called J. F. Colombeau's
generalized numbers (J. F. Colombeau [2], § 2.1). Recently,
M.Oberguggenberger and M. Kunzinger [11] proved that the algebras of
the type ${\cal G}(\Omega)$ can be canonically imbedded in the ring $\bar{\C}^{\os_c}$
of pointwise functions from $\os_c$
into $\bar{\C}$, where $\os_c$ is a particular
set of generalized points such that $\Omega\subset\os_c\subset\bar{\R}^d$.
Thus, the generalized functions in ${\cal G}(\Omega)$ can also
be viewed as pointwise functions.
\vskip 0,2 true cm

   The difference between the algebra of asymptotic functions $\re$
and J. F. Colombeau's algebras of generalized functions of the type
${\cal G}(\Omega)$ can be best expressed by comparing their scalars:
\vskip 0,2 true cm

    {\bf a)} The set of the scalars $\bar{\C}$ of ${\cal G}(\Omega)$ is a 
non-archimedean ring with zero divisors, containing a copy of $\C$. In
contrast, the set of the scalars $^\rho\C$ of $\re$ is an algebraically
closed non-archimedean field, containing a copy of $\C$. In addition,
$^\rho\C$ is Cantor complete in the sense that every sequence of closed
circles in $^\rho\C$ with the finite intersection property has a non-empty
intersection in $^\rho\C$. Finally, $^\rho\C$ contains a canonical copy of 
the T. Levi-Civit\`a field of asymptotic series with complex coefficients
(that property gives the name "asymptotic numbers" for the elements of 
$^\rho\C$).
\vskip 0,2 true cm

    {\bf b)} The set of the real scalars $\bar{\R}$ of ${\cal G}(\Omega)$ is a
partially ordered non-archimedean ring with zero divisors
containing a copy of $\R$. In contrast, the set of the real scalars $^\rho\R$
of $\re$ is a real closed non-archimedean field extension of $\R$
which is Cantor complete in the sense that every sequence of closed
intervals in $^\rho\R$ with the finite intersection property has a non-empty 
intersection in $^\rho\R$. Moreover, $^\rho\R$ contains a canonical copy of
the T. Levi-Civit\`a field of asymptotic series with real
coefficients (A. H. Lightstone and A. Robinson [6], p.93).
\vskip 0,2 true cm

    The improvement of the properties of the scalars $^\rho\R$ and $^\rho\R$ of
$\re$ compared with the properties of the scalars $\bar{\C}$
and $\bar{\R}$ in J. F. Colombeau's theory is the main reason
to involve the methods of nonstandard analysis in our approach.
\vskip 0,5 true cm

\section{Preliminaries: Asymptotic Numbers, Asymptotic Vectors
and Asymptotic Functions}

      Our  framework is a nonstandard model of the real numbers $\R$,
with degree of saturation larger than card$(\N)$. Throughout this
paper $\Omega$ denotes an open subset of $\R^d$. We denote by $^\ast\R$,
$^\ast\R_+$, $^\ast\C$, $^\ast\e$ and $^\ast{\cal D}(\Omega)$
the nonstandard extensions of $\R$, $\R_+$, $\C$, $\e$ and
${\cal D}(\Omega)$, respectively. In particular, $^\ast\e$ consists of pointwise
functions of the form $f$: $^\ast\Omega \to\,^\ast\C$. If $X$ is a set of complex numbers
or a set of (standard) functions, then $^\ast X$ will be its nonstandard
extension and if $f:X\to Y$ is a (standard) mapping, then $^\ast f:\,^\ast X \to\,^\ast Y$
will be its nonstandard extension. For integration in $^\ast\R^d$ we use
the $\ast$-Lebesgue integral. For a very short introduction to
nonstandard analysis we refer to the Appendix in T. Todorov [15].
For a more detailed exposition we recommend T. Lindstr{\o}m [5],
where the reader will find many references to the subject.
\vskip 0,2 true cm

    The  purpose of this section is to present a summary  of  those
definitions and results of {\it Nonstandard Asymptotic Analysis} which we
need for our exposition. For more details we refer to A. Robinson
[13], A.H. Lightstone and A. Robinson [6], W.A.J. Luxemburg [7], Li
Bang-He [4] and V. Pestov [12]. For the algebra of asymptotic
functions and its applications we refer to M. Oberguggenberger [8],
T. Todorov [14], R. F. Hoskins and J. Sousa  Pinto [3], M.
Oberguggenberger and T. Todorov [10], M. Oberguggenberger [9] and
T. Todorov [15].
\vskip 0,2 true cm

{\bf Part 1.} {\it Asymptotic  Numbers:} Let $\rho\in\,^\ast\R$ be a positive
infinitesimal. We shall keep $\rho$ fixed in what follows. Following
A. Robinson [13] we define the field of A. Robinson's complex $\rho$-{\it
asymptotic numbers} as the factor space
$^\rho\C = {\cal M}_\rho(^\ast\C)/{\cal N}_\rho(^\ast\C)$, where
\begin{equation}
  {\cal M}_\rho(^\ast\C)=\{\xi\in\,^\ast\C:|\xi|<\rho^{-n}\quad {\text{for some}}\quad
   n\in\N\},  
\end{equation}
\begin{equation}
  {\cal N}_\rho(^\ast\C)=\{\xi\in\,^\ast\C:|\xi|<\rho^n \quad {\text{for all}}\quad n\in\N\},
\end{equation}
are the sets of the r-moderate and $\rho$-null numbers in $^\ast\C$,
respectively. We define the embedding $\C\subset\,^\rho\C$ by $c\to q(c)$, where
$q: {\cal M}_\rho(^\ast\C)\to\,^\rho\C$ is the quotient mapping.
If ${\cal A}\subseteq\,^\ast\C$, we let
\begin{equation}
    {\cal M}_\rho({\cal A})=\{\xi\in{\cal A}:|\xi|<\rho^{-n}\quad {\text{for some}} \quad n\in\N\},
\end{equation}
and $^\rho{\cal A}=q[{\cal M}_\rho({\cal A})]$. The set $^\rho{\cal A}\subseteq\,^\rho\R$
is called the $\rho$-extension of ${\cal A}$
in $^\rho\R$. In the particular case ${\cal A}=\,^\ast A$, where $A\subseteq\C$, we shall write
simply $^\rho A$ instead of the more precise $^\rho(^\ast A)$. The field of $A$.
{\em Robinson's real $\rho$-asymptotic numbers} is defined by
$^\rho\R = q[{\cal M}_\rho(^\ast\R)]$. The imbedding
$\R\subset\,^\rho\R$ is defined by $r\to q(r)$. The asymptotic number $s=q(\rho)$ is
the {\cal scale} of $^\rho\R$; it is a positive infinitesimal, i.e. $s\in\,^\rho\R,\,s>0$
and $s\approx 0$. Finally, we observe that if $\alpha\in\,^\ast\C$, then
$\alpha\in{\cal M}_\rho(^\ast\C)$ iff $|\alpha|\in{\cal M}_\rho(^\ast\R)$ and,
similarly, $\alpha\in{\cal N}_\rho(^\ast\C)$ iff $|\alpha|\in{\cal N}_\rho(^\ast\R)$.
\vskip 0,2 true cm

    It is easy to check that ${\cal N}_\rho(^\ast\C)$ is a maximal ideal in
${\cal M}_\rho(^\ast\C)$
and hence $^\rho\C$ is a field. Also $^\rho\R$ is a real closed totally ordered
nonarchimedean field (since $^\ast\R$ is a real closed totally ordered
field) containing $\R$ as a totally ordered subfield. Thus, it follows
that $^\rho\C=\,^\rho\R(i)$ is an algebraically closed field, where $i = \sqrt{-1}$.
Finally $^\rho\R$ is Cantor complete in the sense that every sequence of
closed intervals in $^\rho\R$ with the finite intersection property has a
non-empty intersection (since $^\ast\R$ has this property). Unless it is
said otherwise, we shall always topologize $^\rho\R$ and $^\rho\C$ by the order
topology in $^\rho\R$. The sequence of intervals:
\begin{equation}
         I_n = \{x\in\,^\rho\R: -s^n<x<s^n\},
\end{equation}
forms a base for the neighborhoods of the zero in the {\em interval
topology} in  $^\rho\R$.
\vskip 0,2 true cm

{\bf Part 2}. {\it Asymptotic Vectors}: Let $^\ast\R^d$   and  $^\rho\R^d$    be   the
corresponding $d$-dimensional vector spaces over the fields  $^\ast\R$   and
$^\rho\R$, repsectively. We shall write always  $\|\ \|$  for the norms  of
$\R^d$,  $^\ast\R^d$  and $^\rho\R^d$, leaving to the reader to decide from the context
which of these norms we have in mind. Notice that $\|\ \|$ takes non-negative 
values in $\R$, $^\ast\R$ and $^\rho\R$, respectively. In the case  $d=1$
this  norm  reduces to the absolute value in $\R$,  $^\ast\R$ and  $^\rho\R$,
respectively. It is clear that $\R^d$ is a vector subspace of both  $^\ast\R^d$
and  $^\rho\R^d$  over $\R$, in symbols, $\R^d\subset\,^\ast\R^d$ and 
$\R^d\subset\,^\rho\R^d$. We  call  the
elements  in $^\rho\R^d$  {\it asymptotic points}  or {\it asymptotic vectors}.  Unless
it is said otherwise, we shall always topologize $^\rho\R^d$ by the 
{\it product-interval  topology},  i.e. by the product topology in $^\rho\R^d$, 
generated by the order topology in  $^\rho\R$. The sequence of balls:
\begin{equation}
        B_n=\{x\in\,^\rho\R^d: \|x\| < s^n\},\quad n\in \N,
\end{equation}
forms a base for the neighborhoods of the zero in  $^\rho\R^d$.
\vskip 0,2 true cm

    Let   $x\in\,^\rho\R^d$.  Then  $x$  is called {\it infinitesimal,  finite  or
infinitely  large}  if its  $\rho$-norm  $\|x\|$  is infinitesimal,  finite
or  infinitely large in  $^\rho\R$, respectively. We denote by  ${\cal F}(^\ast\R^d)$ and
${\cal F}(^\rho\R^d)$  the  sets  of  the  finite  elements  of   $^\ast\R^d$   and   $^\rho\R^d$,
respectively, and by ${\cal I}(^\ast\R^d)$ and  ${\cal I}(^\rho\R^d)$ we denote the sets  of  the
infinitesimals points in $^\ast\R^d$ and $^\rho\R^d$, respectively. If $\Omega$ is an open
set of $\R^d$, then $\os=\Omega + {\cal I}(^\ast\R^d)$ denotes the set of the
{\it nearstandard  points  of}    $^\ast\Omega$  in $^\ast\R^d$ and
$^\rho\os = \Omega+{\cal I}(^\rho\R^d)$ denotes the set of the {\it nearstandard points of} $^\ast\Omega$
 in $^\rho\R^d$, in symbols,
\begin{equation}
    \os=\Omega + {\cal I}(^\ast \R^d) = \{x + \xi \,: \, x \in \Omega, \; \xi \in {\cal I}(^\ast\R^d)\},
\end{equation}
\begin{equation}
    ^\rho\os=\Omega + {\cal I}(^\ast\R^d) = \{x+h \,: \, x \in \Omega, \; h\in {\cal I}(^\ast\R^d)\},
\end{equation}
where the sum in  $x+\xi$ is in $^\ast\R^d$,  while the sum in  $x+h$ is
in  $^\rho\R^d$.  We observe that $\Omega\subset\os\subseteq\,^\ast\Omega$ and
$\Omega\subset\,\ros\subseteq\,\ro$. In the particular case $\Omega=\R^d$ we have
\begin{equation}
     \R^d+{\cal I}(^\ast\R^d)={\cal F}(^\ast\R^d)\ \text{and}\ \R^d+{\cal I}(^\rho\R^d)={\cal F}(^\rho\R^d).
\end{equation}
We  denote  by  $(^\rho\C)^{^\rho\os}$ or $(^\rho\C)^{\Omega+{\cal I}(^\rho\RR^d)}$ the  
ring (under the pointwise operations) of the functions from $^\rho\os$ into $^\rho\C$.
\vskip 0,2 true cm

    The  purpose of the next several lines is to connect  the  sets
$\os$ and $^\rho\os$. First, we observe  that  the
vector  space $^\rho\R^d$ can, alternatively, be constructed as the  factor
space  $^\rho\R^d = {\cal M}_\rho(^\ast\R^d)/{\cal N}_\rho(^\ast\R^d)$, where
\begin{equation}
     {\cal M}_\rho(^\ast\R^d)=\{\xi\in ^\ast\R^d:\|\xi\| < \rho^{-n}\quad{\text{for some}}\quad n\in\N\},
\end{equation}
\begin{equation}
     {\cal N}_\rho(^\ast\R^d)=\{x \in ^\ast\R^d:\|x\| < \rho^{-n}\quad{\text {for all}}\quad n\in\N\},
\end{equation}
are   the  sets  of  the  $\rho$-moderate  and  $\rho$-null  points  in  $^\ast\R^d$,
respectively.  We denote by  $q^d:{\cal M}_\rho(^\ast\R^d)\to\,^\rho\R^d$  the  corresponding
quotient  mapping. We  let $q^1=q$ in the case $d=1$. The  imbedding
$\R^d\subset\,^\rho\R^d$  is reproduced by  by  $x\to q^d(x)$. If 
${\cal A}\subseteq\,^\ast\R^d$, we let
\begin{equation}
{\cal M}_\rho({\cal A})=\{\xi\in {\cal A}:\| \xi \|<\rho^{-n}\quad {\text {for some}} \quad n\in\N\},
\end{equation}
and $^\rho{\cal A}=q[{\cal M}_\rho({\cal A})]$. The set  $^\rho{\cal A}\subseteq\,
^\rho\R^d$ is called the $\rho$-extension  of
${\cal A}$ in  $^\rho\R^d$. In the particular case ${\cal A}=\,^\ast A$, where 
$A\subseteq\R^d$, we shall write
simply  $^\rho A$  instead  of  the more precise $^\rho(^\ast A)$.  We  observe  that
$\ros$ is the $\rho$-extension of  $\os$   and  also we have:
\begin{equation}
     \xi\in\os \Leftrightarrow q^d(\xi)\in\,\ros.
\end{equation}
For more details about the vector spaces over the field $^\rho\R$ we refer
to W.A.J. Luxemburg [7].
\vskip 0,2 true cm

{\bf Part 3.} {\it Asymptotic Functions}: The algebra of {\it asymptotic
functions}
$\re$ on $\Omega$ is the factor space $\re = {\cal M}_\rho(^\ast\e)/{\cal N}\rho(^\ast\e)$, where
\begin{equation}
  {\cal M}_\rho(^\ast\e)=\{f\in\,^\ast\e:(\forall\alpha\in\N^d_0)(\forall\xi\in\os)(\partial^\alpha
  f(\xi)\in{\cal M}_\rho(^\ast\C))\},
\end{equation}
\begin{equation}
{\cal N}_\rho(^\ast\e)=\{f\in\,^\ast\e:(\forall\alpha\in\N^d_0)(\forall\xi\in\os)(\partial^\alpha
  f(\xi)\in{\cal N}_\rho(^\ast\C))\}.
\end{equation}
The  functions  in  ${\cal M}_\rho(^\ast\e)$ are called  $\rho$-moderate  (or, simply,
{\it moderate)}  and those in ${\cal N}_\rho(^\ast\e)$ are called $\rho${\it -null functions} (or,
simply,  {\it null functions}).  The {\it canonical  (differential   ring)
imbedding} $\e\subset\,\re$ is defined by $f\to\sigma(f)$, where $\sigma(f)=Q_\Omega(^\ast f)$
and $Q_\Omega:{\cal M}_\rho(^\ast\e)\to\,\re$ is  the quotient  mapping.  For  the
construction of an imbedding of the Schwartz distributions:
$\e\subset{\cal D'}(\Omega)\subset\,\re$, discussed in the Introduction of this  paper,  we
refer to (M. Oberguggenberger and T. Todorov [10]).
\vskip 0,5 true cm

\section{Asymptotic Functions as Pointwise Functions}

   We define an evaluation of the asymptotic functions in $\re$ and
show  that  this evaluation coincides with the usual evaluation  in
the  particular  case of classical smooth functions.  Moreover,  we
construct  a  ring  imbedding of $\re$ into the ring  of  pointwise
functions  $(^\rho\C)^{\ros}$, thus, showing that the  asymptotic
functions   can   be  characterized  as  pointwise  functions.   In
particular, the Schwartz distributions also can be characterized by
their values. In what follows we shall use the notations in Section
1.
\vskip 0,2 true cm

    Recall  that  every asymptotic function is, by  definition,  an
equivalence  class  of  nonstandard  internal  functions  in  $^\ast\e$
(Section 1, Part 3).  On the other hand, the functions in $^\ast\e$ are
pointwise  functions from $^\ast\Omega$ into $^\ast\C$. We should  notice,  however,
that  the asymptotic functions $Q_\Omega(f)$ can not inherit literally  the
values of its representatives $f$ because the factorization in  $Q_\Omega(f)$
destroys this evaluation.
\vskip 0,2 true cm

{\bf(2.1) Definition} {\it (Values)}: Let  $F=Q(f)$ be an asymptotic function
in $\re$  with a representative $f\in{\cal M}_\rho(^\ast\e)$. Let $x = q(\xi)$ be an
asymptotic  point  in $\ros$ with a  representative $\xi\in\os$.
The asymptotic number  $q(f(\xi))\in\,^\rho\C$  is called the
{\it value}  of  $F$  at  $x$, in symbols,
\[
    F(x) =  q(f(\xi))\quad {\text or} \quad Q(f)(q(\xi)) = q(f(\xi)).
\]
We  should mention that the evaluation of $f$ at $\xi$ is in  the
framework of  $^\ast\C$.
\vskip 0,2 true cm

{\bf (2.2) Lemma}{\it (Correctness)}:  Let
$\alpha,\,\beta\in\os,\; \|\alpha-\beta\|\in{\cal N}_\rho(^\ast\R)$ and
$f,\,g\in{\cal M}_\rho(^\ast\e)$, $f - g\in{\cal N}_\rho(^\ast\e)$.
Then $f(\alpha)-g(\beta)\in{\cal N}_\rho(^\ast\C)$ (or, equivalently,
$|f(\alpha) - g(\beta) | \in{\cal N}_\rho(^\ast\R))$.
\vskip 0,2 true cm

{\bf Proof}: By the Mean Value Theorem, applied by Transfer Principle (T.
Todorov [15], p. 686), there exists  $\gamma\in\,^\ast\R^d$  between $\alpha$ and
$\beta$  (in the  sense  that
$\gamma = \alpha + \eps(\beta-\alpha)$ for some $\eps\in\,^\ast\R, 0\le\epsilon\le 1)$
 such that
\[
   |f(\alpha)-f(\beta)|=|\nabla f(\gamma)\cdot(\alpha-\beta)|\le\|\nabla
   f(\gamma)\|\,\|\alpha-\beta\|\,.
\]
The  point $\gamma$ is in $\os$, since $\alpha\approx\beta$, by assumption,
thus, $\|\nabla\, f(\gamma)\|$ is a moderate number in ${\cal M}_\rho(^\ast\R)$, by the definition
of   ${\cal M}_\rho(^\ast\e)$. It follows that the right hand side of the above inequality  is
in ${\cal N}_\rho(^\ast\R)$, since ${\cal N}_\rho(^\ast\R)$ is an ideal in
${\cal M}_\rho(^\ast\R)$. On  the other hand, we have
\[
   \|f(\alpha)-g(\beta)\,\|=\|\, f(\alpha)-f(\beta)+f(\beta)-g(\beta)\,\| \le\,
   \|f(\alpha)-f(\beta)\| +
\]
\[
 + \|f(\beta)-g(\beta)\|\,\le\, \|\nabla\, f(\gamma)\|\,\|\alpha
   -\beta\|+\|f(\beta)-g(\beta)\| \in {\cal N}_\rho(^\ast\R),
\]
as required, since $\| f(\beta) - g(\beta)\|\in{\cal N}_\rho(^\ast\R)$, by assumption.
\hfill $\Box$
\vskip 0,2 true cm

    Recall  that the algebra of standard smooth functions  $\e$  is
imbedded in $\re$ in a canonical way, in symbols, $\e\subset\,\re$, by
the mapping $f \to \sigma(f)$, where  $\sigma(f) = Q(^\ast f)$. Recall as well that
$\C\subset\,^\rho\C$
by the mapping  $c\to q(c)$. In what follows we shall identify $c$
with its image  $q(c)$, writing simply   $c = q(c)$.
\vskip 0,2 true cm

    The  next  result  shows that the concept  of  {\it value}   for  the
asymptotic functions, introduced above, is a generalization of  the
usual evaluation of the classical functions.
\vskip 0,2 true cm

{\bf (2.3) Lemma} {\it (Standard Functions at Standard Points)}:
If  $f\in\e$ and  $x \in\Omega$, then  $\sigma(f)(x) = f(x)$.
\vskip 0,5 true cm

{\bf Proof:} $\sigma(f)(x)= Q(^\ast f)(q(x)) = q(^\ast f(x)) = q(f(x)) = f(x)$   since
$^\ast f$   is  an  extension  of  $f$  and  $f(x)$ is  a  standard  (complex)
number.
\hfill$\Box$
\vskip 0,2 true cm

The next two lemmas are in the framework of  $^\ast\R^d$  and $^\ast \C$.
\vskip 0,2 true cm

{\bf (2.4) Lemma:} Let $f\in {\cal M}_\rho(^\ast\e)$ and  $\xi\in\os$ and
$\chi\in{\cal N}_\rho(^\ast\R^d)$.  Then:\\
\phantom{i}(i) $\;\| f(\xi+\chi) - f(\xi)\| \in {\cal N}_\rho(^\ast\C)$.
\vskip 0,3 true cm

          (ii) $\dfrac{\|f(\xi+\chi)-f(\xi)-\nabla f(\xi)\cdot\chi\|}
                      {\|\chi\|}\;   \in{\cal N}_\rho(^\ast\C),\quad\chi\ne 0.$
\vskip 0,2 true cm

{\bf Proof:} (i) By the Mean Value Theorem, applied by Transfer Principle
(T.  Todorov [15], p. 686) there exists  $\eps\in\,^\ast\R,\; 0\le\eps\le 1$,  such
that
\[
 \|f(\xi+\chi)-f(\xi)\|  =  \|\nabla f(\xi+\eps \chi)\cdot \chi \|
\]
It  follows   $\|f(\xi+\chi)-f(\xi)\|\le\|\nabla\,f(\xi+\eps\chi)\|\;\|\chi\|$.
Notice that  $\xi+\eps\chi$   is also in  $\os$ (since $\eps\chi$ is  an
infinitesimal point) which implies $\|\nabla—f(\xi+\eps\chi)\|\in{\cal M}_\rho(^\ast\R)$ since
$f\in{\cal M}_\rho(^\ast\e)$,  by  assumption. Now, the result  follows  since
${\cal N}_\rho(^\ast\R)$ is an ideal in ${\cal M}_\rho(^\ast\R)$.
\vskip 0,2 true cm

    (ii)  By the Taylor Theorem, applied by Transfer Principle,
there exists  $\eps\in\,^\ast\R,\;0\le\eps\le 1$,  such that
\[
  \frac{\|f(\xi+\chi)-f(\xi)-\nabla f(\xi)\cdot\chi\|}
       {\|\chi\|}
       = \frac 12 \left\|\sum_{|\alpha|=2}\,\partial\,^\alpha
       f(\xi+\eps\chi) \frac{\chi^\alpha}{\|\chi\|}\,\right\|
\]
As  before,  $\|\partial ^\alpha f(\xi+\eps\chi)\|\in{\cal M}_\rho(^\ast\R)$ and also we have
$\chi^\alpha/\|\chi\|\in{\cal N}_\rho(^\ast\R)$ due to  $|\alpha |=2.$
Thus, it follows that the right hand side of the above equality is in 
${\cal N}_\rho(^\ast\R)$, since ${\cal N}_\rho(^\ast\R)$ is an ideal in 
${\cal M}_\rho(^\ast\R)$.
\hfill$\Box$
\vskip 0,1 true cm

    In what follows  $K\subset\subset\Omega$ means that  $K$  is a compact subset of
$\Omega$.
\vskip 0,3 true cm

{\bf(2.5) Lemma:}  Let $f\in{\cal M}_\rho(^\ast\e)$.  Then:\\
{\phantom{i}} (i)  $(\forall K\subset\subset\Omega)(\forall n\in\N)(\exists m\in\N)
                    (\forall\xi\in\,^\ast K)(\forall\chi\in\,^\ast\R)$
    \vskip 0,1 true cm

{\hskip 0,5 true cm} 
$[(\|\chi\|<\rho^m)\Rightarrow(\| f(\xi+\chi)-f(\xi)\|<\rho^n)].$
\vskip 0,3 true cm
   (ii) $(\forall K\subset\subset\Omega)(\forall n\in\N)(\exists
   m\in\N)(\forall\alpha\in\,^\ast K)(\forall\chi\in\,^\ast\R)$
   \vskip 0,1 true cm
   {\hskip 0,5 true cm}
        $\left[(0<\|\chi\|<\rho^m) \Rightarrow\quad
        \dfrac{\|f(\xi+\chi)-f(\xi)-\nabla f(\xi)\cdot\chi\|}{\|\chi\|}\,<\,\rho^n\right].$
\vskip 0,2 true cm

{\bf Proof:}  (i) Suppose (for contradiction) that (i) fails, i.e.
\[
   (\exists K \subset\subset\Omega)(\exists n\in\N)(\forall m\in\N)(S_m\ne\not0),
   \]
where
\[
    S_m=\left\{(\xi,\chi)\in\,^\ast K\,\times\,^\ast\R^d:(\|\chi\|<\rho^m)\quad{\text
    and}\quad
    (\| f(\xi+\chi)-f(\xi)\|\ge \rho^n)\right\}.
\]
It  is  clear  that $S_m\supset S_{m+1}\supset ...$, hence, by the  Saturation
Principle (T. Todorov [15], p. 687), there exists a pair  $(\xi_0,\chi_0)\in S_m$
for all  $m\in \N$.  The latter means $\|\chi_0\|\in{\cal N}_\rho(^\ast\R)$  and
$\| f(\xi_0+\chi_0)-f(\chi_0)\|\notin{\cal N}_\rho(^\ast\R)$, contradicting the assumption
$f\in{\cal M}_\rho(^\ast\e)$ in view of Lemma 2.4.
\vskip 0,2 true cm

    (ii)   As  above, suppose (for contradiction) that (ii)  fails,
i.e.
$(\exists K\subset\subset\Omega)(\exists n\in\N)(\forall m\in\N)(T_m \neq
\emptyset)$, where
\begin{equation*}
\begin{split}
T_m = \biggl\{&(\xi,\chi)\in\,^\ast K\,\times\,
^\ast\R^d\,:\\ &0<\|\chi\|<\rho^m\texts{and}
\dfrac{\| f(\xi+\chi)-f(\xi)-\nabla f(\xi)\cdot
   \chi\|}{\|\chi\|}\ge\rho^n\biggr\}.
\end{split}
\end{equation*}

As  before,  we  observe that $T_m\subset T_{m+1}....$ ,  hence,  by  the
Saturation Principle, there exists a pair
$(\xi_0 , \chi_0) \in T_m$  for all  $m\in\N$. The latter means  $\|\chi_0\| \in{\cal N}_\rho(^\ast\R)$
and
\[
  \dfrac{\| f(\xi_0 +\chi_0)-f(\xi_0)-\nabla f(\xi_0)\cdot \chi_0\|}{\|\chi_0\|}
  \;\not\in\;{\cal N}_\rho(^\ast\R),
\]
contradicting the assumption  $f\in{\cal M}_\rho(^\ast\e)$ in view of Lemma  2.4.
\hfill$\Box$
\vskip 0,3 true cm

{\bf (2.6) Theorem:}  Let  $F\in\,^\rho{\cal E}(\Omega)$, $\partial ^\alpha F$
 be its partial derivatives  and $\nabla\, F$  be  the  corresponding gradient of
 $F$ in  $^\rho{\cal E}(\Omega)$.  Let  $x\in\,^\rho\os$
 be a finite asymptotic number and  $F(x)$, $(\partial ^\alpha F)(x)$ and
 $(\nabla \, F)(x)$ be the values of $F$,  $\partial ^\alpha F$ and $\nabla\,F$ at the point $x$,
respectively. Then:\\
{\phantom i}(i) $\lim\limits_{^\rho\RRR^d\ni h\to 0}\;F(x+h) = F(x).$
\vskip 0,3 true cm

    (ii) $\lim\limits_{^\rho\RRR^d\ni h\to 0}\;\dfrac{\|F(x+h)-F(x)-(\nabla F)(x)\cdot h\|}{\|h\|} = 0$
\vskip 0, true cm
where  both  limits are in the product-interval  topology  of   $^\rho\R^d$ (1.5).
\vskip 0,3 true cm

{\bf Proof:}  We  have  $F=Q(f), x = q(\xi)$ and $h = q(\chi)$ for some
$f\in{\cal M}_\rho(^\ast{\cal E}(\Omega))$
some $\xi\in\os$ and some $\chi\in{\cal M}_\rho(^\ast\R^d)$. Now, both
(i)  and (ii) follow immediately from the above lemma after  taking
quotient  mappings  $Q_\Omega$ and $q^d$, taking into  account  that  the
sequence of balls  $\{B_n\}$ (1.5) forms a base for the neighborhoods of
the zero in the product-interval topology in   $^\rho\R^d$.
\hfill$\Box$
\vskip 0,2 true cm

{\bf(2.7) Corollary} ({\it The Case} $\Omega= \R$): Let $F \in\,^\rho{\cal E}(\R)$ and $F'$ be its
derivative  in  the  algebra  $^\rho{\cal E}(\R)$. Let  $x\in{\cal F}(^\rho\R)$  be  a  finite
asymptotic number and $F(x)$ and  $F'(x)$  be the values of  $F$  and $F'$
at the point  $x$, respectively. Then:\\
{\phantom i}(i) $\lim\limits_{^\rho\RRR^d\ni h\to 0}\;F(x+h) = F(x).$
\vskip 0,3 true cm

    (ii) $\lim\limits_{^\rho\RRR^d\ni h\to 0}\;\dfrac{F(x+h)-F(x)}{h} = F'(x)$
\vskip 0,3 true cm
where  both  limits are in the interval topology of $^\rho\R$  (1.4)
(or, equivalently, in the valuation-metric topology of  $^\rho\R$ ).
\vskip 0,3 true cm

{\bf Proof:}  The  result  follows as a particular case  from  the  above
theorem for $d = 1$ and $\Omega= R$, taking into account that in this case
we have $\widetilde{\R} = \R+{\cal I}(^\ast\R)={\cal F}(^\ast\R)$  and
$^\rho\widetilde{\R} = \R+{\cal I}(^\rho\R)={\cal F}(^\rho\R)$.
\hfill$\Box$
\vskip 0,3 true cm

{\bf (2.8) Theorem} {\it(The Imbedding)}: Define
$V:\;^\rho{\cal E}(\R) \to (^\rho\C)^{^\rho\os}$,  by $F\to V(F)$, where $V(F)(x)=F(x)$,
$x\in\,\ros$. Then $V$ is a ring imbedding in the sense  that  it
is injective and it preserves the addition and multiplication.
\vskip 0,3 true cm

{\bf Proof:} Suppose that  $V(F) = 0$ in $(^\rho\C)^{\ros}$, i.e. $F(x)=0$
in $^\rho\C$ for all $x\in\,\ros$. We have to show that  $F=0$
in  the algebra $\re$, i.e. to show that  $F=Q_\Omega(f)$ for some $f\in\,^\ast{\cal E}(\R)$
such that  $\partial ^\alpha f(\xi)\in{\cal N}_\rho(^\ast\C)$ for all  $\xi\in\os$  and
all  $\alpha\in\N^d_0$  (1.14)).  By our assumption we have $F(x+h)-F(x)=0$  for all
$x\in\,\ros$ and all $h\in{\cal I}(^\rho\R^d)$
(since $x + h\in\,\ros$ which implies  $F(x+h)=0$  along
with $F(x)=0$). It follows that $\|\nabla F(x)\|=0$  for  all
$x\in\ros$,
by the second part of Theorem 2.6.  In  other
words, $\partial ^\alpha F(x)=0$ for all  $x\in\,\ros$ and all
$\alpha\in\N^d_0$.
Next, observe that each partial  derivative
$\partial ^\alpha F$ with $|\alpha|=1$, also belongs to $\re$, so, we can  repeat  the
above arguments by induction in $|\alpha|$; the result is  $(\partial ^\alpha F)(x)=0$
in $^\rho\C$ for all  $x\in\,\ros$ and all $\alpha\in\N^d_0$.
On  the other hand,  we have  $F=Q_\Omega(f)$ for some $f\in{\cal M}_\rho(^\ast{\cal E}(\R))$ and
$\partial^\alpha F=Q_\Omega(\partial^\alpha f)$. Also, $x\in\,\ros$ means $x=q^d(\xi)$  for
some  $\xi\in\os$  (1.12).  Hence,  we  conclude   that
$q(\partial^\alpha f(\xi))=0$ for  all  $\xi\in\os$ and  all $\alpha\in\N^d_0.$
The latter means $\partial^\alpha f(\xi)\in{\cal N}_\rho(^\ast\C)$ for all $\xi\in\os$
and all $\alpha\in\N^d_0$, as required. To show
the  preservation  of  the ring operations,  let  us  take  another
asymptotic function $G = Q(g)$ with a representative  $g\in{\cal M}_\rho(^\ast{\cal E}(\R)).$
For the addition we have
\[
  (F+G)(x)= Q(f+g)(q(\xi))=q(f(\xi)+g(\xi))=q(f(\xi))+q(g(\xi))=
\]
\[
                = Q(f)(q(\xi)) + Q(g)(q(\xi)) = F(x) + G(x),
\]
and similar for the multiplication.
\hfill$\Box$
\vskip 0,2 true cm

{\bf (2.9)  Corollary} ($C^\infty${\it -Functions)}: Let  $F\in\,\re$ and $V(F)$ be its
image in $(^\rho\C)^{\ros}$. Then $V(F)$  is a pointwise
$C^\infty$-function from  $\ros$ to $^\rho\C$ in the sense that for any  $x\in\,\ros$
and  any  $\alpha\in\N^d_0$ the  partial derivative
$\partial^\alpha\,V(F)$  exists and it is continuous at   $x$   in the product-order
topology in $^\rho\R^d$ (1.5).
\vskip 0,3 true cm

{\bf Proof:}  The result follows directly from Theorem 2.6,  by induction
in $|\alpha|$.
\hfill$\Box$
\vskip 0,2 true cm

\section{Fundamental Theorem of Calculus in $\re$}

      We  prove the Fundamental Theorem of Integration Calculus  in
$\re$ and  show  that the set of the scalars  of  $\re$  coincides
exactly with A. Robinson's field $^\rho\C$ and the set of the real scalars
of $\re$ coincides exactly with A. Robinson's field  $^\rho\R.$
\vskip 0,3 true cm

{\bf{(3.1) Definition}} ({\it Scalars in}  $\re$): We call  $F\in\,\re$ a scalar
in $\re$ if $\partial^\alpha F=0$  in  $\re$ for all $\alpha\in\N^d_0,|a| = 1$,
in short,  if  $\nabla\,F=0$ in  $\re$.
\vskip 0,3 true cm

{\bf (3.2)  Theorem} ({\it{Fundamental Theorem in}}  $\re)$:  Let  $\Omega$  
be an open subset of $\R^d$ which in addition is arcwise connected. 
Let $F\in\re$  be an asymptotic function and  $\nabla\,F$ be  its  
gradient in $\re$. Then the following are equivalent:
\begin{enumerate}
\item[(i)] $\,F$ is a scalar in $\re$, i.e. $\nabla\,F=0$ in $\re$.

\item[(ii)] $(\nabla\,F)(x)=0$ in  $^\rho\C^d$ for all $x\in\,\ros$
(pointwisely) in the sense that
$(\partial^\alpha\,F)(x)=0$ in $^\rho\C$ for all
$x\in\,\re$  and all $\alpha\in\N^d_0$, $|\alpha| = 1$.

\item[(iii)] $F$ is a constant asymptotic function in the sense  that
$F(x) = C$  in $^\rho\C$ for some  $C\in\,^\rho\C$  and all  $x\in\,\ros$.

\item[(iv)] $F=Q_\Omega(c)$ in  $\re$  for some  $c\in{\cal M}_\rho(^\ast\C)$, 
more precisely, $F=Q_\Omega(f_c)$,
defined by $f_c(\xi)=c$ for all $\xi\in\,^\ast\Omega$.
\end{enumerate}
\vskip 0,2 true cm
The connection between (iii) and (iv) is given by  $C = q(c)$,  where
$q:{\cal M}_\rho(^\ast\C)\to\,^\rho\C$ is the quotient mapping (Section 1, Part 1).
\vskip 0,3 true cm

{\bf{Proof:}} (i) $\,\Leftrightarrow\,$(ii) follows immediately from 
Theorem 2.8.
\vskip 0,1 true cm

(ii) $\,\Rightarrow\,$(iii):   Observe   that $V(\|\nabla\,F\|)=0$   in
$(^\rho\C)^{\ros}$ implies $\|\nabla\,F\|=0$ in $\re$, by  Theorem  2.8.
Let $x_1,\,x_2\in\,\ros$. We have to show that $F(x_1)=F(x_2)$
in $^\rho\C$. We need the representatives: we have $F=Q_\Omega(f)$ for some
$f\in{\cal M}_\rho((\e)$ and $x_i=q^d(\xi_i)$,  for some  $\xi_i\in\os,\,i=1$,
2.  Since $\Omega$ is arcwise connected, it follows that  $^\ast\Omega$  is arcwise
connected, by Transfer Principle (T. Todorov [15], p. 686).  Hence,
$\xi_1$ and $\xi_2$ can be connected with a $\ast$-continuous  curve 
$L\subset\,\ros$. It follows
\[
   f(\xi_2)-f(\xi_1) = \int_L\;\nabla\,f(\xi)\cdot d\xi,
\]
by  Transfer Principle. Continuing the arguments, $\|\nabla\,F\|=0$  in
$\re$ implies $\|\nabla\,f\|\in{\cal N}_\rho(^\ast\e)$.  Thus,
\[
  | f(\xi_2)-f(\xi_1)| \le \|\nabla\,f(\xi_0)\|\|\xi_2-\xi_1\|\in{\cal N}_\rho(^\ast\R)
\]
where $\xi_0$  is some point on L. It follows
$f(\xi_2)-f(\xi_1)\in{\cal N}_\rho(^\ast\C)$.
Hence, $F(x_2)-F(x_1)=q(f(\xi_2)-f(\xi_1))=0$, as required.
\vskip 0,2 true cm

      (ii)$\,\Leftarrow\,$(iii):  We have $F(x + h) - F(x) = 0$  for  all 
$x\in\,\ros$ and all $h\in{\cal I}(^\rho\R^d)$, by assumption (since  $x  +  h$
is  also in  $\ros$). It follows $\nabla\,F(x) = 0$ for all  $x\in\,\ros$,
by the (ii) - part of Theorem 2.6.
\vskip 0,2 true cm

    (iii)$\,\Rightarrow\,$(iv): Suppose that  $F = Q_\Omega(f)$ for some
$f\in{\cal M}_\rho(^\ast{\cal E}(\Omega))$
and $C=q(c)$ for some $c\in{\cal M}_\rho(^\ast\C)$. Now, $F(x) = C$ for all 
$x\in\,\ros$  implies
$q(f(\xi)-c)=0$ or, equivalently, $f(\xi)-c\in{\cal N}_\rho(^\ast\C)$, for all  $\xi\in\os$ (1.14).
Hence,  $F-Q_\Omega(f_c)=Q_\Omega(f - f_c) = Q_\Omega(f-c)= 0$,  as required.
\vskip 0,2 true cm

    (iii)$\,\Leftarrow\,$ (iv): We have $F(x) = q(f_c(\xi)) = q(c) = C$, where $x=q(\xi)$.
\hfill$\Box$
\vskip 0,3 true cm

{\bf (3.3) Corollary:}
\vskip 0,2 true cm

   (i) The set of the scalars in $\re$:
   \[
      \{F\in\re:\;\nabla\,F=0\quad {\text in}\quad \re\}
   \]
\vskip 0,1 true cm
is   an   algebraically  closed  non-archimedean  field  which   is
isomorphic to the field of A. Robinson's complex asymptotic numbers
$^\rho\C$ (Section 1, Part 1) under the imbedding ${\cal C}\,:\,F\to F(x)$, where $x$
is  an  arbitrary point in $\ros$. Consequently, $^\rho\C$ is
imbedded in $\re$, in symbols, $^\rho\C\subset\,\re$,  by ${\cal C}\,:\,C\to Q_\Omega(f_c)$,
where $f_c\in\,^\ast\e$ is defined by $f_c(\xi)=c$ for all  $\xi\in\,^\ast\Omega$, where
$C=q(c)$.   The   imbedding   ${\cal C}^{-1}$  preserves   the   addition   and
multiplication.
\vskip 0,2 true cm

   (ii) The set of the real scalars in $\re$:
\[
   \{\pm |F|\,:\,F\in\,\re\quad \text{such that}\quad \nabla\, F=0\quad \text{in}\quad \re \}
\]

is  a  totally  ordered real closed Cantor complete non-archimedean
field  which  is  isomorphic to the field  of  A. Robinson's  real
asymptotic numbers $^\rho\R$ (Section 1, Part 1) under the imbedding
${\cal C}\,:\,\pm|F|\to\pm|F(x)|$, where  $x$ is an arbitrary point in $\ros.$
\vskip 0,2 true cm

{\bf Proof:}  (i) For every scalar $F$ in $\re$ there exists $C\in\,^\rho\C$  such
that   $F(x) = C$  for all $x\in\,\ros$, by Theorem 3.2.  So,
the mapping  $F\to F(x)$ from the constant functions in $\re$ into  $^\rho\C$
is   well   defined.   The  preservation  of   the   addition   and
multiplication follows from Theorem 2.8. Conversely, if  $C = q(c)$
is in $^\rho\C$, then the function $F=Q_\Omega(f_c)$  is the preimage of $C$.
\vskip 0,2 true cm

    (ii)   follows immediately from  (i)  taking into account  that
$^\rho\R=\{\pm|a|:\alpha\in\,^\rho\C\}.$
\hfill$\Box$
\vskip 0,5 true cm

{\bf Acknowledgment:} The author thanks Michael Oberguggenberger for  the
discussion of the preliminary version of this paper.
\vskip 1 true cm


\address{T. D. Todorov}{California Polytechnic State University, Department of
Mathematics, San Luis Obispo, Califirnia 93407, USA}
{ttodorov\char'100polymail.calpoly.edu}

\end{document}